\documentclass[reqno]{amsart}

\usepackage{amssymb}
\usepackage[all]{xy}

\theoremstyle{definition}
\newtheorem{defn}{Definition}[section]
\theoremstyle{remark}
\newtheorem{rem}[defn]{Remark}

\theoremstyle{plain}
\newtheorem{prop}[defn]{Proposition}
\newtheorem{cor}[defn]{Corollary}
\newtheorem{thm}[defn]{Theorem}
\newtheorem{lemma}[defn]{Lemma}

\newcommand{\rank}{\mathrm{rk}}
\newcommand{\Bunstack}{\mathcal{B}un}
\newcommand{\Sympl}{\mathcal{S}ympl}
\newcommand{\Substack}{\mathcal{S}ub}
\newcommand{\Aa}{\mathbb{A}}
\newcommand{\Extstack}{\mathcal{E}xt}
\newcommand{\Gr}{\mathrm{Gr}}
\newcommand{\GL}{\mathrm{GL}}
\newcommand{\SL}{\mathrm{SL}}
\newcommand{\Gp}{\mathrm{Gp}}
\newcommand{\Sp}{\mathrm{Sp}}
\newcommand{\cohom}{\mathrm{H}}
\newcommand{\Hom}{\mathrm{Hom}}
\newcommand{\RHom}{\mathrm{RHom}}
\newcommand{\Ext}{\mathrm{Ext}}
\newcommand{\Hombdl}{\mathit{Hom}}
\newcommand{\Endbdl}{\mathit{End}}
\newcommand{\Aut}{\mathrm{Aut}}
\newcommand{\longto}[1][]{\stackrel{#1}{\longrightarrow}}
\newcommand{\Quot}{\mathrm{Quot}}
\newcommand{\Spec}{\mathrm{Spec}}
\newcommand{\Pic}{\mathrm{Pic}}
\newcommand{\stackM}{\mathcal{M}}
\newcommand{\Symcoarse}{\mathfrak{Sympl}}
\renewcommand{\O}{\mathcal{O}}
\newcommand{\univ}{\mathrm{univ}}
\newcommand{\generic}{\mathrm{generic}}
\DeclareMathOperator{\rad}{rad}
\DeclareMathOperator{\codim}{codim}
\newcommand{\Euniv}{\mathcal{E}^\univ}
\newcommand{\Echeckuniv}{\check{\mathcal{E}}^\univ}
\newcommand{\Etildeuniv}{\widetilde{\mathcal{E}}^\univ}
\newcommand{\Funiv}{\mathcal{F}^\univ}
\newcommand{\Huniv}{\mathcal{H}^\univ}
\newcommand{\stackE}{\mathcal{E}}
\newcommand{\stackF}{\mathcal{F}}
\newcommand{\stackL}{\mathcal{L}}
\newcommand{\stackI}{\mathcal{I}}
\newcommand{\stackH}{\mathcal{H}}
\newcommand{\stackU}{\mathcal{U}}
\newcommand{\stackV}{\mathcal{V}}
\newcommand{\stackW}{\mathcal{W}}
\newcommand{\integers}{\mathbb{Z}}
\newcommand{\dual}{\mathrm{dual}}
\newcommand{\coker}{\mathrm{coker}}
\DeclareMathOperator{\im}{im}
\newcommand{\id}{\mathrm{id}}
\newcommand{\pr}{\mathrm{pr}}
\newcommand{\ad}{\mathrm{ad}}
\newcommand{\Gm}{\mathbb{G}_m}
\newcommand{\BGm}{\mathrm{B}\mathbb{G}_m}
\renewcommand{\P}{\mathrm{P}}
\newcommand{\PP}{\mathbb{P}}

\begin{document}

\title[Moduli of symplectic bundles]{Some moduli stacks of symplectic 
bundles on a curve are rational}

\author[I. Biswas]{Indranil Biswas}

\address{School of Mathematics, Tata Institute of Fundamental
Research, Homi Bhabha Road, Mumbai 400005, India}

\email{indranil@math.tifr.res.in}

\author[N. Hoffmann]{Norbert Hoffmann}

\address{Mathematisches Institut der
Georg-August-Universit\"at, Bunsenstra{\ss}e 3--5,
37073 G\"ottingen\\Germany}

\email{hoffmann@uni-math.gwdg.de}

\subjclass[2000]{14H60, 14A20}

\keywords{Symplectic bundle, moduli stack, rationality}

\date{}

\begin{abstract}
Let $C$ be a smooth projective curve of genus $g \geq 2$ over a field $k$.
Given a line bundle $L$ on $C$, let $\Sympl_{2n, L}$ be the moduli stack
of vector bundles $E$ of rank $2n$ on $C$ endowed with a nowhere 
degenerate
symplectic form $b: E \otimes E \longto L$ up to scalars. We 
prove that
this stack is birational to $\BGm \times \Aa^s$ for some $s$ if
$\deg( E) = n \cdot \deg( L)$ is odd and $C$ admits a rational point $P \in C( k)$
as well as a line bundle $\xi$ of degree $0$ with $\xi^{\otimes 2} \not \cong {\O}_C$.
It follows that the corresponding coarse moduli scheme of Ramanathan-stable
symplectic bundles is rational in this case.
\end{abstract}

\maketitle

\section{Introduction}

Let $C$ be a smooth projective curve of genus $g \geq 2$ over an algebraically closed field $k$.
Soon after the coarse moduli schemes of vector bundles $E$ over $C$ were constructed in the sixties,
the question of their rationality (in the fixed determinant case) was addressed. It is relatively
easy to prove that these moduli schemes are unirational. Newstead proved that if the rank and
degree are mutually coprime, then the moduli scheme is rational under a numerical condition \cite{Ne}.
King and Schofield showed that the assumption that the rank and degree are mutually coprime is enough
to ensure that the moduli scheme is rational \cite{KS}.

This coprime case is also the case where a Poincar\'{e} family of vector bundles
parameterized by the coarse moduli scheme exists. A deeper relation between rationality
and the existence of Poincar\'{e} families appears in the proof of King and Schofield:
Arguing by induction on the rank, they use not only the birational type of the coarse
moduli scheme for lower rank, but also the obstruction class against the existence of
Poincar\'{e} families on it.

We address the same rationality question for moduli spaces of vector bundles with
symplectic structure. Instead of coarse moduli schemes, we prefer to study moduli stacks,
with the aim of showing that they are birational to $\BGm \times \Aa^s$ for some $s$.
The latter means in more classical terms that the corresponding coarse moduli scheme is
rational, and that some open subscheme of it admits a Poincar\'{e} family.

Here are the moduli stacks that we work with: Given a positive integer $n$ and a line bundle
$L$ on $C$, we consider the moduli stack $\Sympl_{2n, L}$ of rank $2n$ vector bundles $E$
on $C$, equipped with a nondegenerate symplectic form $b: E \otimes 
E\longto  L$;
as isomorphisms between two such bundles
$$
(E, b: E \otimes E  \longto  L) ~~\,\,~~\text{~and~}
~~\,\,~~ (E', b': E' \otimes E' \longto  L)
$$
we allow all vector bundle isomorphisms 
between
$E$ and $E'$ that commute with $b$ and $b'$ up to an automorphism of $L$. Such symplectic
bundles $(E, b)$ can also be viewed as principal bundles under a well-known reductive
group $\Gp_{2n}$ which contains the symplectic group $\Sp_{2n}$ as a normal subgroup
with $\Gm$ as the quotient; see Section \ref{bundles} for the details.

Our main result, Theorem \ref{theorem1}, states that the stack $\Sympl_{2n, L}$ is
birational to $\BGm \times \Aa^s$ for some $s$ if $n$ and the degree of $L$ are both odd.
The latter condition ensures the existence of a Poincar\'{e} family on some open subscheme
of the coarse moduli scheme; in that sense, it is analogous to the condition on vector
bundles that their rank and degree be coprime.

For that theorem, we do not have to assume that the ground field $k$
is algebraically closed;
instead, it suffices for us that the curve $C$ admits a rational point 
$P \in C( k)$ and also
a line bundle $\xi$ of degree $0$ with $\xi^{\otimes 2} \not\cong {\O}_C$.

The idea of the proof is to find a \emph{canonical} reduction of the structure group for
every sufficiently general $\Gp_{2n}$-bundle. More precisely, we show that every sufficiently
general symplectic bundle $(E, b)$ admits a canonical line subbundle of $E$. This allows us
to reconstruct $E$ from bundles of lower rank and some appropriate extension data; we manage
to parameterise all these rationally.

\section{Symplectic bundles on a curve} \label{bundles}

Let $k$ be an arbitrary field. Given vector spaces $V$ and $L$ over $k$ with $\dim( L) = 1$,
we call a bilinear form
\begin{equation*}
  b: V \otimes V \longto L
\end{equation*}
\emph{alternating} if $b( v, v) = 0$ for all $v \in V$. This implies that $b( v, w) + b( w, v) = 0$
for all $v, w \in V$. These two conditions are equivalent if the characteristic of $k$ is different from $2$.

The \emph{adjoint} of $b$ is the linear map
\begin{equation*}
  b^{\#}: V \longto \Hom( V, L), \qquad v \mapsto b(\_ \, , v): V 
\longto  L\, .
\end{equation*}
The form $b$ is called \emph{nondegenerate} if $b^{\#}$ is an isomorphism. If $b$ is
nondegenerate, then the alternating form $b$ is called a \emph{symplectic
form}. If $b$ is symplectic, then $\dim( V)$ is finite
and even, say $\dim( V) = 2n$. In this case, the triple $( V, L, b)$ is isomorphic to $k^{2n}$
equipped with the standard symplectic form, so the automorphism group of 
the triple $( V, L, b)$ is isomorphic to the group
\begin{equation*}
  \Gp_{2n}( k) := \{(A, \lambda) \in \GL_{2n}( k) \times k^* \big|
 A^t \left( \begin{smallmatrix} 0 & I_n\\-I_n & 0 \end{smallmatrix} \right) A
  = \lambda \cdot \left( \begin{smallmatrix} 0 & I_n\\-I_n & 0 \end{smallmatrix} \right) \}
\end{equation*}
where $I_n$ denotes the $n \times n$ unit matrix. Varying $k$, we obtain a reductive algebraic group
$\Gp_{2n}$. It occurs in a canonical short exact sequence
\begin{equation*}
  1 \longto \Sp_{2n} \longto \Gp_{2n} \longto[\pr_2] \Gm \longto 1
\end{equation*}
where $\Sp_{2n}$ denotes the usual symplectic group and $\pr_2( A, 
\lambda) := \lambda$. Using the standard
fact that $\Sp_{2n} \subseteq \SL_{2n}$, it is easy to check that
\begin{equation} \label{det}
  \det( A) = \lambda^n \text{ for all } (A, \lambda) \in \Gp_{2n}( k)\, .
\end{equation}

Now let $C$ be a smooth, geometrically irreducible, projective curve of genus $g \geq 2$ over $k$ with a
rational point $P \in C( k)$. We denote by $k_P$ the coherent skyscraper sheaf supported at
$P$ with stalk $k$. A \emph{vector bundle} $E$ on $C$ is a locally free coherent sheaf; a \emph{subbundle} of $E$ is a
coherent subsheaf with torsion-free quotient. Let
$$
\eta_P: \cohom^0( 
E)\longto  E_P
$$
be the canonical evaluation map to the fibre $E_P$ of $E$ at $P$.

We consider (twisted) \emph{symplectic bundles} $E$ on $C$. Here 
`twisted' means that we replace the $1$-dimensional vector space $L$ 
above by a line bundle on $C$; slightly abusing notation, we use
$L$ for denoting this line bundle as well. So a twisted symplectic 
bundle on $C$ consists more precisely of
a rank $2n$ vector bundle $E$ on $C$, a line bundle $L$ on $C$ and a nowhere degenerate symplectic form
\begin{equation*}
  b: E \otimes E \longto L\, .
\end{equation*}
Such symplectic bundles correspond to principal $\Gp_{2n}$-bundles on $C$, as follows: Given a principal
$\Gp_{2n}$-bundle, the associated vector bundles $E$ and $L$ can be constructed by means of the canonical
representations
\begin{equation*}
  \pr_1: \Gp_{2n} \longto \GL_{2n}, \quad (A, \lambda) \mapsto A \quad\text{and}\quad
  \pr_2: \Gp_{2n} \longto \Gm,      \quad (A, \lambda) \mapsto \lambda\, .
\end{equation*}
Conversely, given a rank $2n$ symplectic bundle $(E, b: E \otimes 
E\longto  L)$, the pair $(E, L)$ determines
a principal $(\GL_{2n} \times \Gm)$-bundle, and $b$ determines a reduction of its structure group to $\Gp_{2n}$.
These two constructions are inverses of each other.
Eq. \eqref{det} implies that the symplectic form
$$
b: E \otimes E\longto  L
$$
induces an isomorphism
\begin{equation*}
  b^{\otimes n}: \det( E) \longto[\sim] L^{\otimes n}.
\end{equation*}

Fix a line bundle $L$ on $C$, and fix a positive integer $n$. We denote by
\begin{equation}\label{syla.}
  \Sympl_{2n, L}
\end{equation}
the moduli stack of symplectic bundles $(E, b: E \otimes 
E\longto  L)$ on $C$ of the
above type with $\rank (E)= 2n$; it is more
precisely given
by the following groupoid $\Sympl_{2n, L}( S)$ for every $k$-scheme $S$:
\begin{itemize}
 \item Each object consists of a vector bundle $\stackE$ of rank $2n$
on $C \times_k S$, a line bundle $\stackL$
  on $C \times_k S$ which is locally in $S$ isomorphic to the pullback of $L$ from $C$, and a nowhere
  degenerate symplectic form $b: \stackE \otimes \stackE\longto  
\stackL$.
 \item Each morphism from $(\stackE, b: \stackE \otimes 
\stackE\longto  \stackL)$ to $(\stackE', b': \stackE'
  \otimes \stackE'\longto  \stackL')$ consists of vector bundle 
isomorphisms $\stackE \longto  \stackE'$ and
  $\stackL \longto  \stackL'$ that intertwine $b$ and $b'$.
\end{itemize}
This stack is known to be algebraic (in the sense of Artin) and locally of finite type over $k$. By standard
deformation theory, it is smooth of dimension $n( 2n+1)( g-1) - 1$ over $k$.
We will show that $\Sympl_{2n, L}$ is also irreducible; see Corollary \ref{irreducible}.
We denote by
\begin{equation*}
  \Symcoarse_{2n, L}
\end{equation*}
the corresponding coarse moduli scheme of Ramanathan-stable symplectic bundles $( E, b: E \otimes E \longto L)$ as above.
It is a normal quasi-projective variety of dimension $n( 2n+1)( g-1)$ over $k$; cf. \cite{ramanathan} and \cite{GLSS}.

Every invertible function $f \in \Gamma( S, \O_S^*)$ induces an automorphism of every object
$$
(\stackE, b:
\stackE \otimes \stackE \longto  \stackL)
$$
in $\Sympl_{2n, L}(S)$, acting by multiplication with $f$ on $\stackE$
and by multiplication with $f^2$ on $\stackL$.
This defines a canonical group homomorphism
\begin{equation} \label{scalars}
  \Gamma( S, \O_S^*) \longto \Aut( \stackE, b: \stackE \otimes \stackE 
\longto  \stackL).
\end{equation}
It will be necessary to keep track of these scalar automorphisms systematically.
Some terminology for that purpose is introduced in \cite{par}; for the convenience of the reader,
we repeat the definitions here.
\begin{defn}
  A \emph{$\Gm$-stack} $\stackM = (\stackM,\iota)$ over $k$ consists of an algebraic stack $\stackM$ over $k$
  together with a group homomorphism $\iota( \stackE): \Gamma( S, \O_S^*) \to \Aut_{\stackM( S)}( \stackE)$
  for each $k$-scheme $S$ and each object $\stackE$ of the groupoid $\stackM( S)$ such that the diagrams
  \begin{equation*} \xymatrix{
      \Gamma( S, \O_S^*) \ar[r]^-{\iota( \stackE)} \ar[dr]_{\iota( \stackE')} & \Aut_{\stackM( S)}( \stackE)
      \ar[d]^{\alpha \mapsto \varphi \alpha \varphi^{-1}}\\ & \Aut_{\stackM( S)}( \stackE')
    } \quad\text{and}\quad \xymatrix{
      \Gamma( S, \O_S^*) \ar[r]^-{\iota( \stackE)} \ar[d]_{f^*} & \Aut_{\stackM( S)}( \stackE) \ar[d]^{f^*}\\
      \Gamma( T, \O_T^*) \ar[r]^-{\iota( f^* \stackE)}          & \Aut_{\stackM( T)}( f^* \stackE)
  } \end{equation*}
  commute for each morphism $\varphi:\stackE \to \stackE'$ in $\stackM(S)$ and each $k$-morphism $f:T \to S$.
\end{defn}
For example, the above group homomorphisms \eqref{scalars} turn $\Sympl_{2n, L}$ into a $\Gm$-stack.
\begin{defn}
  A $\Gm$-stack $(\stackM, \iota)$ is a \emph{$\Gm$-gerbe} if $\iota( \stackE): \Gamma( S, \O_S^*) \to \Aut_{
  \stackM(S)}(\stackE)$ is an isomorphism for every $k$-scheme $S$ and every object $\stackE$ of $\stackM( S)$.
\end{defn}
For example, the open substack in $\Sympl_{2n, L}$ of symplectic bundles admitting only scalar automorphisms is a $\Gm$-gerbe.
This open substack is known to be non-empty for all $n \geq 1$.
\begin{defn}
  Let $(\stackM, \iota)$ and $(\stackM', \iota')$ be $\Gm$-stacks over $k$. A $1$-morphism $\Phi: \stackM \to
  \stackM'$ \emph{has weight $w \in \integers$} if the diagram
  \begin{equation*} \xymatrix{
    \Gamma(S,\O_S^*) \ar[r]^-{\iota(\stackE)}\ar[d]_{(\_)^w} &\Aut_{\stackM(S)}(\stackE)\ar[d]^{\Phi(S)}\\
    \Gamma(S,\O_S^*) \ar[r]^-{\iota'(\Phi(\stackE))} &\Aut_{\stackM'( S)}( \Phi( \stackE))
  } \end{equation*}
  commutes for every $k$-scheme $S$ and every object $\stackE$ of the groupoid $\stackM( S)$.
\end{defn}

\section{Reduction of structure group}

\begin{lemma} \label{lin_alg}
  Let $L$ and $I \subseteq H \subseteq V$ be finite dimensional $k$-vector spaces with $\dim(L) = \dim(I) =1$
  and $\dim( H) = \dim( V) - 1$. Then there is a canonical bijection between the following two collections:
  \begin{itemize}
   \item nondegenerate symplectic forms $b: V \otimes V \longto  
L$ such that $I^{\perp} = H$, and
   \item nondegenerate symplectic forms $c: H/I \otimes H/I 
\longto  L$ together with an isomorphism of
    short exact sequences
    \begin{equation*} \xymatrix{
      0 \ar[r] & H/I \ar[d]^{c^{\#}} \ar[r] & V/I \ar[d]^f \ar[r] & V/H \ar[d]^{\bar{f}} \ar[r] & 0\\
      0 \ar[r] & \Hom(H/I, L)        \ar[r] & \Hom(H, L)   \ar[r] & \Hom( I, L)          \ar[r] & 0
    } \end{equation*}
    in which $c^{\#}$ is the adjoint of $c$.
  \end{itemize}
\end{lemma}
\begin{proof}
  One direction is easy: Given the form $b$, its restriction to $H \otimes H$ induces a nondegenerate
  symplectic form $c$ on $H/I$, and the required isomorphism $(c^{\#}, f, \bar{f})$ of short exact sequences
  is induced by the adjoint $b^{\#}: V \longto  \Hom( V, L)$ of 
$b$.

  For the converse direction, assume that a pairing $c$ and an isomorphism of short exact
  sequences $(c^{\#}, f, \bar{f})$ as above are given. It is easy to
  see that the diagram of canonical linear maps
  \begin{equation*} \xymatrix{
    H \otimes H \ar[r] \ar[d] & \Lambda^2 H \ar[d]\\V \otimes H \ar[r] & \Lambda^2 V
  } \end{equation*}
  is cocartesian; hence the linear maps $\Lambda^2 H \longto  L$ 
and $V \otimes H \longto  L$ given by $c$ and $f$ are
  both induced by a unique linear map $\Lambda^2 V \longto  L$. 
The latter defines a symplectic form $b$ with the
  required properties.
\end{proof}

\begin{defn} \label{substack}
  For a given line bundle $I$ on $C$, we denote by
  \begin{equation*}
    \Substack_{1, I}( \Euniv) \longto[\Phi_I] \Sympl_{2n, L}
  \end{equation*}
  the moduli stack of rank $2n$ symplectic bundles $(E, b: E \otimes E 
\longto  L)$ (see \eqref{syla.}), together with a vector
  subbundle of $E$ isomorphic to $I$.
\end{defn}
More precisely, this moduli stack $\Substack_{1, I}( \Euniv)$ is given by the following groupoid
$\Substack_{1, I}( \Euniv)(S)$ for each $k$-scheme $S$:
\begin{itemize}
 \item Its objects consist of an object $(\stackE, b: \stackE \otimes 
\stackE \longto  \stackL)$ in $\Sympl_{2n, L
  }(S)$ and a subbundle $\stackI \subseteq \stackE$ that is locally in $S$ isomorphic to the pullback of $I$.
 \item Morphisms from $(\stackE, b: \stackE \otimes \stackE 
\longto  \stackL, \stackI \subseteq \stackE)$ to
  $(\stackE', b': \stackE' \otimes \stackE' \longto  \stackL', 
\stackI' \subseteq \stackE')$ consist of three
  vector bundle isomorphisms
$$\stackE \longto  \stackE'\, , ~~\,\,~~
\stackL \longto  \stackL' ~~\,\,~~\text{~and~} ~~\,\,~~
\stackI \longto  \stackI'
$$
that commute with $b$ and $b'$ and with the inclusions $\stackI 
\subseteq \stackE$ and $\stackI' \subseteq \stackE'$.
\end{itemize}
Forgetting the subbundle $\stackI$ defines the canonical $1$-morphism $\Phi_I$ above.

\begin{lemma}
  $\Substack_{1, I}( \Euniv)$ is an algebraic stack locally of finite type over $k$.
\end{lemma}

\begin{proof}
  The $1$-morphism $\Phi_I$ is of finite type and representable, namely by appropriate locally closed
  subschemes of relative $\Quot$-schemes. The lemma thus follows from the corresponding statement about
  $\Sympl_{2n, L}$.
\end{proof}
We consider $\Substack_{1, I}( \Euniv)$ as a $\Gm$-stack in the canonical way that makes $\Phi_I$ a
morphism of weight $1$.

\begin{defn} \label{def:extstack}
  For a line bundle $I$ on $C$, we denote by
  \begin{equation} \label{extstack}
    \Extstack( \Funiv, I) \longto \Sympl_{2n-2, L}
  \end{equation}
  the moduli stack of rank $2n-2$ symplectic bundles $(F, c: F \otimes F 
\longto  L)$, together with a vector
  bundle extension of $F$ by $I$
  \begin{equation} \label{IHF}
    0 \longto I \longto H \longto F \longto 0\, .
  \end{equation}
\end{defn}
More precisely, this moduli stack $\Extstack( \Funiv, I)$ is given by the following groupoid
$\Extstack( \Funiv, I)( S)$ for each $k$-scheme $S$:
\begin{itemize}
 \item Its objects consist of an object $(\stackF, c: \stackF \otimes 
\stackF \longto  \stackL)$ in
  $\Sympl_{2n-2, L}( S)$ together with an exact sequence of vector bundles on $C \times_k S$
  \begin{equation*}
    0 \longto \stackI \longto \stackH \longto \stackF \longto 0
  \end{equation*}
  such that $\stackI$ is locally in $S$ isomorphic to the pullback of $I$.
 \item Its morphisms consist of four vector bundle isomorphisms
$$
\stackF 
\longto  \stackF'\, ,\, \stackL \longto \stackL'\, ,
\, \stackI \longto  \stackI'\,~ \text{~and~}\,~ 
\stackH \longto  \stackH'
$$
that commute with all the given maps.
\end{itemize}
Forgetting the extension by $\stackI$ defines the canonical $1$-morphism in \eqref{extstack}.

\begin{lemma} \label{extlemma}
  (i) $\Extstack( \Funiv, I)$ is an algebraic stack locally of finite type over $k$.

  (ii) The morphism \eqref{extstack} in Definition \ref{def:extstack} is smooth, surjective, and all its fibres are irreducible.
\end{lemma}

\begin{proof}
  (i) Put $L' := I \otimes L^{\otimes (n-1)}$, and let $\Bunstack_{2n-1, 
L'}$ denote the moduli stack of vector bundles $H$
of rank $2n-1$
over $C$ with $\det( H) \cong L'$. We have a canonical $1$-morphism
  \begin{equation*}
    \Extstack( \Funiv, I) \longto \Bunstack_{2n-1, L'}
  \end{equation*}
  that sends the exact sequence in \eqref{IHF} to the vector bundle $H$. This morphism is of finite type and
  representable. It is
represented by appropriate locally closed subschemes of 
iterated relative $\Quot$-schemes that
  parameterise the quotients $H \twoheadrightarrow F$ and $\Lambda^2 F \twoheadrightarrow L$. Since
  $\Bunstack_{2n-1, L'}$ is known to be algebraic and locally of finite type over $k$, the same follows
  for $\Extstack( \Funiv, I)$.

  (ii) It is a direct consequence of \cite[Lemma 1.10]{par}.
\end{proof}
We consider $\Extstack( \Funiv, I)$ as a $\Gm$-stack in the canonical way that makes \eqref{extstack} a
morphism of weight $1$.

\begin{defn}
  For a line bundle $I$ on $C$, the canonical $1$-morphism $\Psi_I$
  \begin{equation*} \xymatrix{
    \Substack_{1, I}( \Euniv) \ar[d] \ar[r]^{\Psi_I} & \Extstack( \Funiv, I) \ar[d] \\
    \Sympl_{2n, L} & \Sympl_{2n-2, L}
  } \end{equation*}
  sends each triple $(\stackE, b: \stackE \otimes \stackE 
\longto  
\stackL, \stackI \subseteq \stackE)$ to the
  vector bundle $\stackF := \stackI^{\perp}/\stackI$, equipped with the symplectic form $c: \stackF \otimes
  \stackF \longto  \stackL$ induced by $b$, and the vector 
bundle extension
  \begin{equation*}
    0 \longto \stackI \longto \stackH := \stackI^{\perp} \longto \stackF \longto 0. 
  \end{equation*}
\end{defn}
Our next aim is to prove that $\Psi_I$ is smooth and surjective with 
irreducible fibres. For that purpose, we
relate it to the following stacks of lifted vector bundle extensions.

\begin{defn}
  Let $A$ be a finitely generated $k$-algebra. Given an exact sequence 
  \begin{equation*}
    0 \longto K \longto E \longto Q \longto 0
  \end{equation*}
  and a morphism $f: \widetilde{K} \longto  K$ of vector bundles 
on $C \times_k \Spec( A)$, we denote by
  \begin{equation*}
    \Extstack( Q, \widetilde{K} \longto[f] K)/E
  \end{equation*}
  the moduli stack over $A$ of lifted vector bundle extensions
  \begin{equation*} \xymatrix{
    0 \ar[r] & \widetilde{K} \ar[r] \ar[d]^f & \widetilde{E} \ar[r] \ar[d] & Q \ar[r] \ar@{=}[d] & 0\\
    0 \ar[r] &        K  \ar[r]          &        E  \ar[r]        & Q \ar[r]            & 0.
  } \end{equation*}
\end{defn}
More precisely, $\Extstack( Q, \widetilde{K} \longto[f] K)/E( S)$ denotes the following groupoid for each
$A$-scheme $\pi: S \longto  \Spec( A)$:
\begin{itemize}
 \item Its objects consist of a vector bundle $\widetilde{\stackE}$ on $C \times_k S$ and a commutative
  diagram with exact rows
  \begin{equation*} \xymatrix{
    0 \ar[r] &\pi^* \widetilde{K} \ar[r] \ar[d]_{\pi^* f} &\widetilde{\stackE} \ar[r] \ar[d]
             &\pi^* Q \ar[r] \ar@{=}[d] &0\\
    0 \ar[r] &\pi^*        K  \ar[r]          &\pi^*        E  \ar[r]        &\pi^* Q \ar[r] &0.
  } \end{equation*}
 \item Its morphisms are the vector bundle isomorphisms 
$\widetilde{\stackE} \longto  \widetilde{\stackE'}$ that
  commute with all the given maps.
\end{itemize}

By definition, we have the cartesian square of canonical $1$-morphisms
\begin{equation} \label{Extlift} \xymatrix{
  \Extstack( Q, \widetilde{K} \longto[f] K)/E \ar[r] \ar[d] & \Extstack( Q, \widetilde{K})
\ar[d]^{f_*}\\
  \Spec( A) \ar[r]^{c_E} & \Extstack( Q, K) 
} \end{equation}
in which $\Extstack( Q, \widetilde{K})$ and $\Extstack( Q, K)$ are moduli stacks over $A$ of vector bundle
extensions (cf. \cite[Example 1.9]{par}), and $c_E$ is the classifying morphism of the given extension $E$.
Using \cite[Lemma 1.10]{par}, this implies in particular that the stack 
$\Extstack(Q, \widetilde{K} \longto  K)/E$
is algebraic and that the stack $\Extstack(Q, \widetilde{K} 
\longto  K)/E$
is of finite type over $k$.

\begin{lemma}\label{lemma2.8}
  Suppose that the homomorphism of vector bundles
$$
f: \widetilde{K} \longto  K
$$
is surjective. Then the structural morphism
  \begin{equation*}
    \Extstack( Q, \widetilde{K} \longto[f] K)/E \longto \Spec( A)
  \end{equation*}
  is smooth and surjective with all fibres irreducible.
\end{lemma}

\begin{proof}
  It suffices to show that the $1$-morphism $f_*$ in \eqref{Extlift} is smooth and surjective with
  irreducible fibres. Using \cite{ega3}, we can represent the morphism
  \begin{equation*}
    f_*: \RHom( Q, K) \longto \RHom( Q, \widetilde{K})
  \end{equation*}
  in the derived category of finitely generated $A$-modules by a chain morphism
  \begin{equation*} \xymatrix{
    V^0 \ar[d]_{\delta} \ar[r]^{f^0} & \widetilde{V}^0 \ar[d]^{\widetilde{\delta}}\\
    V^1                 \ar[r]_{f^1} & \widetilde{V}^1
  } \end{equation*}
  of length one complexes $V^{\bullet}$ and $\widetilde{V}^{\bullet}$ that consist of vector bundles $V^0$,
  $V^1$ and $\widetilde{V}^0$, $\widetilde{V}^1$ on $\Spec(A)$. The proof of \cite[Lemma 1.10]{par} yields
  $1$-isomorphisms
  \begin{equation*}
    V^1/V^0 \longto[\sim] \Extstack( Q, K) \qquad\text{and}\qquad
      \widetilde{V}^1/\widetilde{V}^0 \longto[\sim] \Extstack( Q, \widetilde{K})
  \end{equation*}
  where $V^1/V^0$ and $\widetilde{V}^1/\widetilde{V}^0$ are the Picard stacks over $\Spec( A)$
  associated to $V^{\bullet}$ and $\widetilde{V}^{\bullet}$ \cite[Exp. XVIII, 1.4]{sga4};
  by their construction, the diagram
  \begin{equation*} \xymatrix{
    V^1/V^0 \ar[r]^-{\sim} \ar[d]_{(f^1, f^0)} & \Extstack( Q, K) \ar[d]^{f_*}\\
    \widetilde{V}^1/\widetilde{V}^0 \ar[r]^-{\sim} & \Extstack( Q, \widetilde{K})
  } \end{equation*}
  commutes up to a $2$-isomorphism. Therefore, it remains to show that $(f^1, f^0)$ is smooth and surjective
  with irreducible fibres.

  The diagram of Picard stacks
  \begin{equation*} \xymatrix{
    (V^1 \oplus \widetilde{V}^0) \big/ V^0 \ar[r] \ar[d]_{(f^1, -\widetilde{\delta})}
    & V^1/V^0 \ar[d]^{(f^1, f^0)}\\ \widetilde{V}^1 \ar[r] & \widetilde{V}^1/\widetilde{V}^0
  } \end{equation*}
  is easily checked to be cartesian. Since the canonical morphisms
  \begin{equation*}
    V^1 \oplus \widetilde{V}^0 \longto (V^1 \oplus \widetilde{V}^0) \big/ V^0
      \qquad\text{and}\qquad \widetilde{V}^1 \longto \widetilde{V}^1/\widetilde{V}^0
  \end{equation*}
  are smooth and surjective with irreducible fibres, it suffices to prove the same for
  \begin{equation*}
    (f^1, -\widetilde{\delta}): V^1 \oplus \widetilde{V}^0 \longto \widetilde{V}^1.
  \end{equation*}
  This is now simply a morphism of vector bundles, so we just have to 
show that it is surjective. For that we first note that the
Nakayama's lemma allows us to assume that the ground ring $A$ is a 
field, say our base field $k$. In this
  case, the cokernel of $(f^1, -\widetilde{\delta})$, by the choices of $V^{\bullet}$ and
  $\widetilde{V}^{\bullet}$, is isomorphic to the cokernel of the $k$-linear map
  \begin{equation*}
    f_*: \Ext^1_{\O_C}( Q, K) \longto \Ext^1_{\O_C}( Q, \widetilde{K})\, ;
  \end{equation*}
  hence this cokernel vanishes if $f$ is surjective because $C$ is a smooth curve.
\end{proof}

\begin{cor} \label{sub2extlemma}
  The canonical $1$-morphism $\Psi_I$ is smooth, surjective, and all its fibres are irreducible.
\end{cor}
\begin{proof}
  Let $A$ be a finitely generated $k$-algebra, and let
  \begin{equation*}
    c_{\stackH}: \Spec( A) \longto  \Extstack( \Funiv, I)
  \end{equation*}
  be the classifying morphism of an object
  \begin{equation} \label{IHF_Ext}
    (\stackF, \qquad c: \stackF \otimes \stackF \longto \stackL,
      \qquad 0 \longto \stackI \longto \stackH \longto[p] \stackF \longto 0)
  \end{equation}
  in $\Extstack( \Funiv, I)( \Spec( A))$. From Lemma \ref{lin_alg} we know that the objects in
  $\Substack_{1, I}( \Euniv)( \Spec( A))$ over \eqref{IHF_Ext} correspond to lifted vector bundle extensions
  \begin{equation*} \xymatrix{
    0 \ar[r] & \stackH \ar[r] \ar[d]^{c^{\#} \circ p} & \stackE \ar[r] \ar[d]
             & \stackI^{\dual} \otimes \stackL \ar[r] \ar@{=}[d] & 0\\
    0 \ar[r] & \stackF^{\dual} \otimes \stackL \ar[r] & \stackH^{\dual} \otimes \stackL \ar[r]
             & \stackI^{\dual} \otimes \stackL \ar[r] & 0;
  } \end{equation*}
  furthermore, the analogous statement holds after any base change $\pi: 
S \longto  \Spec( A)$. This means that
  the diagram
  \begin{equation*} \xymatrix{
    \Extstack \big(
      \stackI^{\dual} \otimes \stackL, \stackH \xrightarrow{c^{\#} \circ p} \stackF^{\dual} \otimes L
    \big) \big/ \stackH^{\dual} \otimes \stackL \ar[r] \ar[d] & \Substack_{1, I}( \Euniv) \ar[d]^{\Psi_I}\\
    S \ar[r]^{c_{\stackH}} & \Extstack( \Funiv, I)
  } \end{equation*}
  is cartesian. Thus the corollary follows from Lemma \ref{lemma2.8}.
\end{proof}

\begin{cor}
  The stacks $\Extstack( \Funiv, I)$ and $\Substack_{1, I}( \Euniv)$ are both smooth.
\end{cor}
\begin{proof}
 This follows from the combination of Lemma \ref{extlemma}(ii), 
Corollary \ref{sub2extlemma}, and the
  smoothness of the stack $\Sympl_{2n-2, L}$.
\end{proof}

\begin{cor} \label{irreducible}
  The stacks $\Sympl_{2n, L}$, $\Substack_{1, I}( \Euniv)$ and $\Extstack( \Funiv, I)$ are all irreducible.
  In particular, they are all non-empty.
\end{cor}
\begin{proof}
  We argue by induction on $n$. For $n = 0$, $\Sympl_{0, L} \cong \Spec( k)$ is irreducible.

  For the induction step, let us assume $n \geq 1$ and that $\Sympl_{2n-2, L}$ is non-empty and irreducible.
  According to Lemma \ref{extlemma} and Corollary \ref{sub2extlemma}, first $\Extstack( \Funiv, I)$ over
  $\Sympl_{2n-2, L}$ and then $\Substack_{1, I}( \Euniv)$ over $\Sympl_{2n, L}$ are non-empty and irreducible
  as well, for every line bundle $I$. Hence $\Sympl_{2n, L}$ is non-empty.

  Let $[E_1]$ and $[E_2]$ be an arbitrary pair of points in $\Sympl_{2n, L}$ corresponding to the symplectic
  bundles $E_1$ and $E_2$, respectively. It is easy to see that both $E_1$ and $E_2$ have a line subbundle
  isomorphic to $I$ if $\deg( I) \ll 0$ is sufficiently negative. In this case, both points $[E_1]$ and
  $[E_2]$ are in the image of the irreducible stack $\Substack_{1, I}( \Euniv)$ and hence they are in the
  same component of $\Sympl_{2n, L}$. This proves that $\Sympl_{2n, L}$ 
is indeed irreducible. Therefore, the proof is complete by induction.

\end{proof}

\section{An auxiliary rationality result}

Let $\stackM$ be an irreducible algebraic stack over $k$, endowed with the structure of a $\Gm$-stack.
We say that $\stackM$ is \emph{rational} as a $\Gm$-stack if it is
birational to $\BGm \times \Aa^s$ for some $s$. Here birational means that the two stacks contain non-empty
open substacks which are $1$-isomorphic. We say that $\stackM$ is \emph{unirational} if it admits a dominant
$1$-morphism from a dense open subscheme of $\Aa^s$ for some $s$.

For every vector bundle $\stackV$ on $\stackM$, we denote by $\P \stackV = \Gr_1( \stackV)$ the projective
bundle of lines in the fibres of $\stackV$, and by $\PP \stackV = \Gr_{\rank(\stackV)-1}(\stackV)$ the
projective bundle of hyperplanes in the fibres of $\stackV$ (cf. \cite[Section 4]{par}). The birational
type of such projective bundles will in general depend on the action of the scalar automorphisms in
$\stackM$ on the fibres of $\stackV$. This action is encoded in the notion of weight \cite[p. 526]{KS};
for the convenience of the reader, we repeat here the general definition \cite[Definition 2.6]{par}.
\begin{defn}
  Let $( \stackM, \iota)$ be a $\Gm$-stack over $k$. A vector bundle $\stackV$ on $\stackM$ \emph{has weight
  $w \in \integers$} if the diagram
  \begin{equation*} \xymatrix{
    \Gamma(S,\O_S^*) \ar[r]^-{\iota(\stackE)}\ar[d]_{(\_)^w} &\Aut_{\stackM(S)}(\stackE)\ar[d]^{\stackV(S)}\\
    \Gamma(S,\O_S^*) \ar[r]^-{\cdot \id_{\stackV( \stackE)}} &\Aut_{\O_S}( \stackV( \stackE))
  } \end{equation*}
  commutes for every $k$-scheme $S$ and every object $\stackE$ of the groupoid $\stackM( S)$.
\end{defn}
If the $\Gm$-stack $\stackM$ is rational and the vector bundle $\stackV$ has some weight $w \in \integers$
in this sense, then both projective bundles $\P \stackV$ and $\PP \stackV$ are
again rational $\Gm$-stacks; this can be seen as follows:

We may assume $\stackM = \BGm \times \Aa^s$ without loss of generality. Pulling back the tautological line
bundle of weight $1$ on $\BGm$, we obtain a line bundle of weight $1$ on $\stackM$. Since tensoring $\stackV$
with a line bundle does not change $\P \stackV$ or $\PP \stackV$, this reduces us to the case $w
= 0$. Then
the vector bundle $\stackV$ and hence also the projective bundles $\P \stackV$ and $\PP \stackV$ are trivial
over some open substack $\emptyset \neq \stackU \subseteq \stackM$ due to \cite[Corollary 3.8]{par}. This
implies that $\P \stackV$ and $\PP \stackV$ are indeed both rational as $\Gm$-stacks.

\begin{lemma} \label{nonspecial_extension}
  Given line bundles $L_1$ and $L_2$ on $C$, there is an extension
  \begin{equation} \label{L_2byL_1}
    0 \longto L_1 \longto E \longto L_2 \longto 0
  \end{equation}
  such that the connecting homomorphism $\delta$ in its long exact cohomology sequence
  \begin{equation*}
    0 \longto  \cohom^0( L_1) \longto  \cohom^0( E) 
\longto  
\cohom^0( L_2) \longto[\delta]
          \cohom^1( L_1) \longto  \cohom^1( E) \longto  
\cohom^1( 
L_2) \longto  0
  \end{equation*}
  has maximal rank, more precisely $\rank( \delta) = \min \{\dim \cohom^0( L_2), \dim \cohom^1( L_1)\}$.
\end{lemma}
\begin{proof}
  Serre duality allows us to assume
  \begin{equation} \label{leq}
    \dim \cohom^0(L_2) \leq \dim \cohom^1(L_1)
  \end{equation}
  without loss of generality; we then have to show that $\delta$ is injective.

  All such extensions \eqref{L_2byL_1} are classified by the affine space $\Ext^1( L_2, L_1)$; those with
  injective connecting homomorphism $\delta$ form an open subscheme $U \subseteq \Ext^1( L_2, L_1)$. We will
  prove $U( k) \neq \emptyset$ by estimating the dimension of the complement.

  If $\delta$ is not injective, then some section $s \in \cohom^0( L_2)$ can be lifted to $E$; this means
  that the extension class $[E] \in \Ext^1( L_2, L_1)$ is annihilated by the homomorphism
  \begin{equation} \label{s*}
    s^*: \Ext^1( L_2, L_1) \longto \Ext^1( \O, L_1) \cong \cohom^1( L_1).
  \end{equation}
  Now $s^*$ is surjective because its cokernel embeds into $\Ext^2( L_2/\O, L_1) = 0$, so
  \begin{equation} \label{codim}
    \codim( \ker (s^*) \subseteq \Ext^1( L_2, L_1)) = \dim \cohom^1( L_1).
  \end{equation}
  Such sections $s$ are parameterized by the projective space $\P \cohom^0( L_2)$, whose dimension
  is smaller than the codimension in \eqref{codim} since we have assumed \eqref{leq}. Therefore,
  \begin{equation*}
    \dim( \Ext^1( L_2, L_1) \setminus U) < \dim \Ext^1( L_2, L_1)
  \end{equation*}
  and consequently $U \neq \emptyset$. If $k$ is infinite, then the non-empty open subscheme $U$ of an affine
  space automatically contains a $k$-rational point, and we are done.

  So suppose that $k$ is finite with $q$ elements. Since \eqref{s*} is a surjective $k$-linear map, the
  cardinalities of these vector spaces then satisfy
   \begin{equation*}
     \frac{\# \ker(s^*)}{\# \Ext^1(L_2,L_1)} = \frac{1}{q^d} \qquad\text{with}\qquad d := \dim \cohom^1(L_1).
  \end{equation*}
  On the other hand, the number of such sections $s$ up to $k^*$ is
  \begin{equation*}
    \# \P \cohom^0( L_2) \leq \# \P \cohom^1( L_1) = \frac{q^d - 1}{q - 1} < q^d.
  \end{equation*}
  Hence $\Ext^1(L_2, L_1) \setminus U$ contains less $k$-rational points than $\Ext^1( L_2, L_1)$; consequently,
  $U( k) \neq \emptyset$ holds for finite fields $k$ as well.
\end{proof}

\begin{cor} \label{nonspecial_new}
  Let $n \geq 1$ be given, and let $L$ and $I$ be line bundles on $C$.
  \begin{itemize}
   \item[(i)] If $2 \deg( I) < \deg( L) + 2 - 2g$, then there is a rank $2n$ symplectic bundle
    \begin{equation*}
      (E, b: E \otimes E \longto  L)
    \end{equation*}
    on $C$ which admits a vector subbundle isomorphic to $I$ and satisfies
    \begin{equation*}
      \Ext^1_{\O_C}( I, E) = 0.
    \end{equation*}
   \item[(ii)] If $2 \deg( I) < \deg( L) + 2g - 2$, then there is a rank $2n$ symplectic bundle
    \begin{equation*}
      (F, c: F \otimes F \longto  L)
    \end{equation*}
    on $C$ which satisfies
    \begin{equation*}
      \Hom_{\O_C}(F, I) = 0 \qquad\text{and} \qquad \Ext^1_{\O_C}( F, I) \neq 0
    \end{equation*}
    as well as
    \begin{equation*}
      \Hom_{  \O_C}( I^{\dual} \otimes L, F) = 0 \qquad\text{and}\qquad
      \Ext^1_{\O_C}( I^{\dual} \otimes L, F) \neq 0.
    \end{equation*}
  \end{itemize}
\end{cor}
\begin{proof}
  It suffices to treat the special case of $n = 1$, since the general case immediately follows from it by
  taking the fibrewise orthogonal direct sum of $n$ copies.

  (i) By assumption, the line bundle $I' := I^{\otimes -2} \otimes L$ has degree $\geq 2g - 1$, so $\cohom^1(
  I') = 0$ by Clifford's theorem, and also we have $\dim \cohom^0( I') 
\geq g = \dim \cohom^1( \O)$ by Riemann-Roch. Thus,
according to Lemma \ref{nonspecial_extension}, there is a vector bundle 
extension
  \begin{equation*}
    0 \longto \O \longto E' \longto I' \longto 0
  \end{equation*}
whose connecting homomorphism $\delta: \cohom^0( I') \longto  
\cohom^1( 
\O)$ is surjective; this implies that
  $\cohom^1( E') = \cohom^1( I') = 0$. Now $E := E' \otimes I$, equipped with the symplectic form $b$
  given by $\det( E) \cong I' \otimes I^{\otimes 2} \cong L$, has the required properties.

  (ii) The two $\Ext$-groups in question are nonzero for every rank $2$ vector bundle $F$ with $\det(F)\cong L
  $ due to Riemann-Roch. The two $\Hom$-groups in question are 
isomorphic via $c^{\#}:F \longto  F^{\dual} \otimes
  L$, so it suffices to construct one such symplectic bundle $F$ with $\Hom( I^{\dual} \otimes L, F) = 0$.

  By assumption, the line bundle $L' := I^{\otimes 2} \otimes L^{\dual} \otimes \omega_C^{\dual}$ has
  negative degree, so $\cohom^0( L') = 0$, and $\dim \cohom^1( L') \geq g = \dim \cohom^0( \omega_C)$ by
  Riemann-Roch. Thus, according to Lemma \ref{nonspecial_extension}, 
there is a vector bundle extension
  \begin{equation*}
    0 \longto L' \longto F' \longto \omega_C \longto 0
  \end{equation*}
  whose connecting homomorphism $\delta: \cohom^0( 
\omega_C)\longto  \cohom^1( L')$ is injective; this implies
  $\cohom^0( F') = \cohom^0( L') = 0$. Now $F := F' \otimes I^{\dual} \otimes L$, equipped with the
  symplectic form $c$ given by $\det( F) \cong L' \otimes \omega_C \otimes (I^{\dual} \otimes L)^{\otimes 2}
  \cong L$, has the required properties.
\end{proof}

\begin{prop} \label{unirat}
  Let $n \geq 1$ be given, and let $L$ be a line bundle over our curve $C$.
  \begin{itemize}
   \item[(i)] If $I$ is a line bundle on $C$ with $2 \deg( I) < \deg( L) + 2g - 2$, then the algebraic stack
    $\Extstack( \Funiv, I)$ over $\Sympl_{2n-2, L}$ is rational as a $\Gm$-stack.
   \item[(ii)] If $I$ is a line bundle on $C$ with $2 \deg( I) < \deg( L)$, then the algebraic stack
    $\Substack_{1, I}( \Euniv)$ over $\Sympl_{2n, L}$ is rational as a $\Gm$-stack.
   \item[(iii)] If $\stackV$ is a vector bundle on some non-empty open substack $\stackU \subseteq \Sympl_{2n,
    L}$ with odd weight $w$ and rank
    \begin{equation*}
      \rank( \stackV) \geq
        \begin{cases} n & \text{ if $\deg( L)$ is odd,}\\2n & \text{ if $\deg( L)$ is even,} \end{cases}
    \end{equation*}
    then the projective bundle $\P \stackV$ over $\Sympl_{2n, L}$ is rational as a $\Gm$-stack.
  \end{itemize}
\end{prop}
\begin{proof}
  We first show that (i) implies (ii) and (ii) implies (iii) for any fixed $n \geq 1$.

  (i) $\Rightarrow$ (ii): We consider the strictly commutative diagram of $1$-morphisms
  \begin{equation*} \xymatrix{
    \Substack_{1, I}( \Euniv) \ar[d]_{\Psi_I} \ar@{}[r]|-{\displaystyle \supseteq}
      & \Psi_I^{-1}( \stackU_1) \ar[r]^-{\Gamma_1} \ar[d]_{\Psi_I}
      & \P \stackW \subseteq \P\Ext( I^{\dual} \otimes L, \Huniv) \ar[dl]\\
    \Extstack( \Funiv, I) \ar@{}[r]|-{\displaystyle \supseteq} & \stackU_1 &
  } \end{equation*}
  in which
  \begin{itemize}
   \item $\stackU_1 \subseteq \Extstack( \Funiv, I)$ is the open substack of all rank $2n-2$ symplectic
    bundles $(F, c: F \otimes F \longto  L)$ together with a 
\emph{nontrivial} vector bundle extension
    \begin{equation} \label{IHF_generic}
      0 \longto I \longto H \longto[p] F \longto 0
    \end{equation}
    with the property that $\Hom_{\O_C}(I^{\dual} \otimes L, H) = 0$,
   \item $\Ext( I^{\dual} \otimes L, \Huniv)$ denotes the vector bundle of weight $1$ on $\stackU_1$ whose
    fibre over such an extension \eqref{IHF_generic} is the vector space
    \begin{equation*}
      \Ext^1_{\O_C}( I^{\dual} \otimes L, H)\, ,
    \end{equation*}
    the isomorphisms in $\Extstack( \Funiv, I)$ acting only on the second variable $H$,
   \item $\stackW \subseteq \Ext( I^{\dual} \otimes L, \Huniv)$ is the subbundle whose fibre over such an
    extension \eqref{IHF_generic} is the inverse image, under the canonical surjection
    \begin{equation*}
     (c^{\#} \circ p)_*: \Ext^1_{\O_C}( I^{\dual} \otimes L, H) \longto
                         \Ext^1_{\O_C}( I^{\dual} \otimes L, F^{\dual} \otimes L)\, , 
    \end{equation*}
    of the line spanned by the extension class of $H^{\dual} \otimes L$, and
   \item $\Gamma_1$ is defined by sending each triple
    \begin{equation*}
      E, \qquad b: E \otimes E \longto L, \qquad I \subseteq E
    \end{equation*}
    to the class of the extension
    \begin{equation*}
      0 \longto H := I^{\perp} \longto E \longto[b^{\#}] I^{\dual} \otimes L \longto 0.
    \end{equation*}
  \end{itemize}

  It is a straightforward consequence of Lemma \ref{lin_alg} that 
$\Gamma_1$ is a $1$-isomorphism
  onto $\P \stackW$. Corollary \ref{nonspecial_new}(ii) asserts that $\stackU_1 \neq \emptyset$; hence
  $\Psi_I^{-1}( \stackU_1) \neq \emptyset$ according to Corollary \ref{sub2extlemma}. This shows that
  $\Substack_{1, I}( \Euniv)$ is birational to $\P \stackW$ over $\stackU_1$. Since we assume (i), the
  $\Gm$-stack $\stackU_1$ is rational here, so $\P \stackW$ is rational as well; thus (ii) follows.

  (ii) $\Rightarrow$ (iii): Assigning to each object $(\stackE, b: 
\stackE \otimes \stackE \longto  \stackL)$ of
  $\Sympl_{2n,L}(S)$ the restriction of $\stackL$ to the point $P \in C( k)$ defines a line bundle of weight
  $2$ over $\Sympl_{2n,L}$. Since tensoring the given vector bundle $\stackV$ with a line bundle does not
  change $\P \stackV$, it suffices to consider one particular odd weight $w$ in the proposition, say $w=1$.

  We choose a line bundle $I$ on $C$ with degree
  \begin{equation} \label{deg_I}
    \deg( I) = \begin{cases}
      (\deg( L) + 1)/2 - g & \text{ if $\deg( L)$ is odd,}\\
       \deg( L)/2      - g & \text{ if $\deg( L)$ is even,}
    \end{cases}
  \end{equation}
  and consider the strictly commutative diagram of $1$-morphisms
  \begin{equation*} \xymatrix{
    \Substack_{1, I}( \Euniv) \ar[d]_{\Phi_I} \ar@{}[r]|-{\displaystyle \supseteq}
      & \Phi_I^{-1}( \stackU_2) \ar[r]^-{\Gamma_2} \ar[d]_{\Phi_I} & \P\Hom( I, \Euniv) \ar[dl]\\
    \Sympl_{2n, L} \ar@{}[r]|-{\displaystyle \supseteq} & \stackU_2 &
  } \end{equation*}
  in which
  \begin{itemize}
   \item $\stackU_2 \subseteq \Sympl_{2n, L}$ is the open locus of all bundles $E$ with
    $\Ext^1_{\O_C}( I, E) = 0$,
   \item $\Hom( I, \Euniv)$ denotes the vector bundle of weight $1$ on $\stackU_2$ whose fibre over such a
    symplectic bundle $E$ is the vector space $\Hom_{\O_C}( I, E)$,
   \item $\Gamma_2$ sends each object $(\stackE, b, \stackI \subseteq \stackE)$ to the sheaf of all morphisms
    from the pullback of $I$ to $\stackE$ that factor through the subbundle $\stackI \subseteq \stackE$.
  \end{itemize}

  It is straightforward to verify that $\Gamma_2$ is an open immersion, 
more precisely a $1$-isomorphism onto
  the open locus of all nonzero morphisms $I \longto  E$ up to 
$k^*$ whose cokernel is torsion-free. Corollary
  \ref{nonspecial_new}(i) asserts that $\Phi_I^{-1}( \stackU_2) \neq \emptyset$; hence the stack
  $\Substack_{1, I}( \Euniv)$ is birational to the stack $\P\Hom( I, \Euniv)$. Since we assume (ii), the
  former is rational as a $\Gm$-stack; thus the latter is so as well.

  Our assumption on $\rank( \stackV)$ together with our choice in \eqref{deg_I} and the Riemann-Roch
  theorem ensure
  \begin{equation*}
    \rank( \Hom( I, \Euniv)) \leq \rank( \stackV)\, .
  \end{equation*}
  Since both $\Hom( I, \Euniv)$ and $\stackV$ are weight $1$ vector bundles on open substacks of
  $\Sympl_{2n, L}$, the former is a direct summand of the latter on some possibly smaller open substack of
  $\Sympl_{2n, L}$ due to \cite[Lemma 3.10(iv) and Lemma 3.10(v)]{par}. The rationality of
  $\P\Hom( I, \Euniv)$ thus implies the rationality of $\P \stackV$ according to \cite[Lemma 4.5(i)]{par}.
  This shows that (ii) indeed implies (iii).

  Now we can prove the proposition by induction on $n$. For $n=1$, we have $\Sympl_{2n-2, L} \cong \Spec(k)$,
  and over it $\Extstack( \Funiv, I) \cong \BGm$, so (i) holds trivially. For the induction step, we
  consider the strictly commutative diagram of $1$-morphisms
  \begin{equation*} \xymatrix{
    \Extstack( \Funiv, I) \ar[d] \ar@{}[r]|-{\displaystyle \supseteq}
      & \stackM \ar[r]^-{\Gamma_3} \ar[d] & \P\Ext( \Funiv, I) \ar[dl]\\
    \Sympl_{2n-2, L} \ar@{}[r]|-{\displaystyle \supseteq} & \stackU_3 &
  } \end{equation*}
  in which
  \begin{itemize}
   \item $\stackU_3 \subseteq \Sympl_{2n-2, L}$ is the open substack of all $F$ with $\Hom_{\O_C}(F, I) = 0$,
   \item $\Ext( \Funiv, I)$ denotes the vector bundle of weight $-1$ on $\stackU_3$ whose fibre over such a
    symplectic bundle $F$ is the vector space $\Ext^1_{\O_C}( F, I)$,
   \item $\stackM \subseteq \Extstack( \Funiv, I)$ is the open substack in the inverse image of
    $\stackU_3$ where the extension $0 \longto  I 
\longto  H \longto  F 
\longto  0$ does not split, and
   \item $\Gamma_3$ sends every such nonsplit extension to the extension class of $H$.
  \end{itemize}

  Here it is straightforward to verify that $\Gamma_3$ is a 
$1$-isomorphism. Corollary
  \ref{nonspecial_new}(ii) asserts $\stackM \neq \emptyset$; hence $\Extstack( \Funiv, I)$ and $\P\Ext(
  \Funiv, I)$ over $\Sympl_{2n-2, L}$ are birational. By part (iii) of the induction hypothesis, the latter
  is rational as a $\Gm$-stack; hence the former is so as well, which suffices to complete the induction.
\end{proof}

\begin{cor}
  $\Sympl_{2n, L}$ is unirational.
\end{cor}
\begin{proof}
  It is easy to see that the canonical $1$-morphism
  \begin{equation*}
    \Phi_I: \Substack_{1, I}( \Euniv) \longto \Sympl_{2n, L}
  \end{equation*}
  is dominant for every line bundle $I$ on $C$ of sufficiently negative degree $\deg( I) \ll 0$. According
  to part (ii) of Proposition \ref{unirat}, the stack $\Substack_{1, I}( \Euniv)$ is in particular unirational;
  hence $\Sympl_{2n, L}$ is so as well.
\end{proof}

\section{Rationality of the moduli stack}

The main result is proved in this section. We start with a lemma.

\begin{lemma} \label{generic}
  Suppose that $k = \bar{k}$ is algebraically closed. Let $n \geq 2$ be even, and let $L$ be a line bundle
  on $C$ with
  \begin{equation*}
    \deg( L) = 2g-1 \qquad\text{and}\qquad \cohom^1( C, L( -P)) = 0
  \end{equation*}
  for our chosen point $P \in C( k)$. Then there is a
symplectic bundle of rank $2n$
\begin{equation*}
    (E, b: E \otimes E \longto  L)
  \end{equation*}
  on $C$ with the following properties:
  \begin{itemize}
   \item[(i)] $\dim \cohom^0( E) = n$ and $\cohom^1( E) = 0$.
   \item[(ii)] The canonical morphism of vector bundles on $C$
    \begin{equation} \label{nonspecial_map}
      \O_C \otimes_k \cohom^0( E) \longto E
    \end{equation}
    is injective with torsion-free cokernel.
   \item[(iii)] The induced alternating pairing on the $n$-dimensional vector space $\cohom^0( E)$
    \begin{equation} \label{nonspecial_pairing}
      \cohom^0( E) \otimes \cohom^0( E)
        \xrightarrow{\eta_P \otimes \eta_P} E_P \otimes E_P \longto[b_P] L_P
    \end{equation}
    is nondegenerate.
  \end{itemize}
\end{lemma}
\begin{proof}
  It suffices to treat the special case of $n = 2$, since the general case immediately follows from it by
  taking the fibrewise orthogonal direct sum of $n/2$ copies.

  The Riemann-Roch theorem implies $\dim \cohom^0( L) = g = \dim 
\cohom^1( \O_C)$. Thus, according to Lemma
  \ref{nonspecial_extension}, there is a vector bundle extension
  \begin{equation*}
    0 \longto \O_C \longto[i] F \longto[p] L \longto 0
  \end{equation*}
  with $\cohom^0( F) = \cohom^0( \O_C) \cong k$ and $\cohom^1( F) = \cohom^1( L) = 0$. The alternating form
  \begin{equation*}
    c := \det: F \otimes F \longto L
  \end{equation*}
  turns $F$ into a symplectic bundle on $C$. The fibrewise orthogonal direct sum
  \begin{equation*}
    E := F \perp F
  \end{equation*}
  satisfies (i) and (ii) for $n = 2$, but does not satisfy (iii).

  Let $\stackU \subseteq \Sympl_{4, L}$ be the open substack defined by (i) and (ii); let $\stackU' \subseteq
  \stackU$ be the open substack where (iii) holds as well. We just saw $\stackU \neq \emptyset$.

  Suppose $\stackU' = \emptyset$. For every point $[E]$ in $\stackU$, the alternating pairing
  \eqref{nonspecial_pairing} on $\cohom^0( E) \cong k^2$ is then degenerate; in other words, this pairing on
  $k^2$ is zero. Since $\stackU$ is reduced, it follows that \eqref{nonspecial_pairing} vanishes identically
  on $\stackU$. In particular, it vanishes for every infinitesimal 
deformation $\stackE$ of $E = F \perp F$,
meaning for every object $\stackE$ of $\stackU( k[\epsilon])$ with 
$\epsilon^2 = 0$ that satisfy the condition that its reduction modulo
  $\epsilon$ is isomorphic to this $E$. But we construct below a deformation $\stackE$ of $E = F \perp F$ for
  which the pairing \eqref{nonspecial_pairing} does not vanish identically. This contradiction will show
  $\stackU' \neq \emptyset$, proving the lemma.

  Let $\check{F} \subseteq F$ be the inverse image of the coherent subsheaf $L( -P) \subseteq L$; then
  \begin{equation*}
    0 \longto \O_C \longto[i] \check{F} \longto[\check{p}] L( -P) \longto 0
  \end{equation*}
  is exact as well. In particular, $\cohom^0( \check{F}) \subseteq \cohom^0( F)$ contains the nonzero section
  of $F$; hence $\cohom^0( \check{F}) = \cohom^0( F)$. Since $\deg( \check{F}) = \deg( F) - 1$, Riemann-Roch
  implies that $\dim \cohom^1( \check{F}) = \dim \cohom^1( F) + 1 = 1$.

  We claim that the composed map
  \begin{equation} \label{cup}
    \cohom^0( F) \otimes \Ext^1( F, \O_C) \longto[\cup] \cohom^1( \O_C) \longto [i_*] \cohom^1( \check{F})
  \end{equation}
  is nonzero. In fact, the second map $i_*$ is surjective since $\cohom^1( L( -P)) = 0$ by hypothesis,
  and the first map $\cup$ is also surjective because of the commutative diagram
  \begin{equation*} \xymatrix{
    \cohom^0( \O_C) \otimes \Ext^1(  F, \O_C) \ar[d]^{i_* \otimes \id} \ar[r]^{\id \otimes i^*} & 
    \cohom^0( \O_C) \otimes \Ext^1( \O_C, \O_C) \ar[d]^{\cup}\\
    \cohom^0(  F) \otimes \Ext^1(  F, \O_C) \ar[r]^-{\cup} & \cohom^1( \O_C)
  } \end{equation*}
  in which $i_*$ is an isomorphism by the choice of $F$. The other 
vertical map $\cup$ is obviously an
  isomorphism as well, and $i^*$ is surjective due to $\Ext^2( L, \O_C) = 0$.

  Having shown that \eqref{cup} is really nonzero, we can choose a 
class
$$
\alpha \in \Ext^1( F, \O_C)
$$
such that the composed map
  \begin{equation} \label{cup_alpha}
    k \cong \cohom^0( F) \xrightarrow{\_ \cup \alpha}
            \cohom^1( \O_C) \longto[i_*]
            \cohom^1( \check{F}) \cong k
  \end{equation}
  is nonzero. We will use this class $\alpha$ to construct the required deformation $\stackE$ of $E$.

  Let $G \subseteq \Gp_4$ be the closed subgroup given by
  \begin{equation*}
    G( k) := \{ (A, \lambda) \in \Gp_4( k)\mid A \cdot e_1 = e_1 \}\, ,
  \end{equation*}
  where $e_1 \in k^4$ is the first standard basis vector. The trivial line subbundle
  \begin{equation} \label{subline}
    \O_C \cong 0 \oplus \O_C \subseteq F \oplus F = E
  \end{equation}
  defines a reduction of structure group to $G$ for our $\Gp_4$-bundle $E = F \perp F$. Let
  \begin{equation*}
    \ad_G( E) \subseteq \ad_{Gp_4}( E) \subseteq \Endbdl( F \oplus F)
  \end{equation*}
  be the corresponding adjoint bundles. We consider the bundle morphism
  \begin{equation} \label{bundlemorphism}
    \Hombdl( F, \O_C) \longto \Endbdl( F \oplus F)
  \end{equation}
  which sends a local morphism $f: F \longto  \O_C$ to the local 
endomorphism
  \begin{equation} \label{section}
    \begin{pmatrix} 0 & f^t \circ p\\i \circ f & 0 \end{pmatrix}: F \oplus F \longto F \oplus F\, ,
  \end{equation}
  where the local morphism $f^t: L \longto  F$ is the adjoint of 
$f$, defined by the formula $c( v, f^t( \ell))
  = f(v) \cdot \ell$ for $\ell \in L$ and $v \in F$. It is easy to check that \eqref{section} is indeed a
  local section of the subbundle $\ad_G( E) \subseteq \Endbdl( F \oplus F)$, so \eqref{bundlemorphism}
  restricts to a bundle morphism
$$
\Hombdl( F, \O_C) \longto  \ad_G( E)
$$
and hence it induces a linear map
  \begin{equation*}
    \Ext^1( F, \O_C) \longto \cohom^1( C, \ad_G( E))\, .
  \end{equation*}

  By standard deformation theory, the image of $\alpha$ in $\cohom^1( C, \ad_G( E))$ corresponds to an
  infinitesimal deformation $\stackE$ of $E$ as a $G$-bundle and hence in particular as a $\Gp_4$-bundle.
  It remains to check that the deformed symplectic pairing \eqref{nonspecial_pairing} on the free
  $k[ \varepsilon]$-module $\cohom^0( \stackE)$ of rank $2$ is nonzero.

  Let $\check{E} := F \oplus \check{F} \subseteq F \oplus F = E$ be the subsheaf of all sections of $E$ which
  are orthogonal at $P$ to the distinguished trivial line subbundle \eqref{subline}. Our construction of
  $\stackE$ as a deformed $G$-bundle implies that $\stackE$ still contains a distinguished trivial line
  subbundle. So taking all sections orthogonal at $P$ to that defines a subsheaf $\check{\stackE} \subseteq
  \stackE$ which is a deformation of $\check{E} \subseteq E$. We have to show the inequality
$$
\cohom^0( \check{
  \stackE}) \neq \cohom^0( \stackE)\, ,
$$
for it implies that there is a section of $\stackE$ which is at $P$ not
  orthogonal to our distinguished section, so the symplectic pairing in question is nonzero.

  Since $\check{\stackE}$ is a deformation of $\check{E}$, we have an exact sequence of $\O_C$-modules
  \begin{equation*}
    0 \longto \check{E} \longto[\cdot \epsilon] \check{\stackE} \longto \check{E} \longto 0
  \end{equation*}
  and consequently a long exact sequence of $k$-vector spaces
  \begin{equation*}
    0 \longto \cohom^0( \check{E}) \longto \cohom^0( \check{\stackE}) \longto \cohom^0( \check{E})
      \longto[\delta] \cohom^1( \check{E}) \longto \cdots\, .
  \end{equation*}
  By construction, one of the four matrix entries of $\delta$ is the map \eqref{cup_alpha}, so $\delta
  \neq 0$ and
  \begin{equation*}
    \dim_k \cohom^0( \check{ \stackE}) < \dim_k \cohom^0( \check{E}) + \dim_k \cohom^0( \check{E})
      = 2 + 2 = 4\, .
  \end{equation*}
  On the other hand, we have $\dim_k \cohom^0( \stackE) = 2 \dim_k 
\cohom^0( E) = 4$. This follows for example from the fact that
  $\cohom^1( E) = 0$. This proves that $\cohom^0( \check{ \stackE}) \neq 
\cohom^0( \stackE)$, as required.
\end{proof}

\begin{thm}\label{theorem1}
  Let $n \geq 1$ be odd, and let $L$ be a line bundle of odd degree over our curve $C$. Suppose
  that there is a line bundle $I$ on $C$ with
  \begin{equation*}
    \deg( L) + 2 \deg( I) = 2g - 1 \qquad\text{and}\qquad L \otimes I^{\otimes 2} \not\cong \omega_C( P)
  \end{equation*}
  for some rational point $P \in C( k)$. Then the $\Gm$-stack $\Sympl_{2n, L}$ is rational.
\end{thm}
\begin{proof}
  Sending each symplectic bundle $(E, b: E \otimes E \longto  
L)$ to the symplectic bundle
  \begin{equation*}
    (E \otimes I, \qquad b \otimes \id: (E \otimes I) \otimes (E \otimes I) \longto L \otimes I^{\otimes 2})
  \end{equation*}
  we construct a $1$-isomorphism $\Sympl_{2n, L} \longto[\sim] \Sympl_{2n, L \otimes I^{\otimes 2}}$. This
  reduces us to the case $I = \O_C$; in other words, we may assume without loss of generality that
  \begin{equation*}
    \deg( L) = 2g - 1 \qquad\text{and}\qquad L \not\cong \omega_C( P)\, .
  \end{equation*}
  According to Clifford's theorem, these imply that $\cohom^1( C, L( 
-P)) = 0$.

  We consider the moduli stack
  \begin{equation*}
    \check{\Sympl}_{2n, L, P}
  \end{equation*}
  of vector bundles $\check{E}$ on $C$ of rank $2n$
endowed with an alternating bilinear form
  \begin{equation*}
    \check{b}: \check{E} \otimes \check{E} \longto L
  \end{equation*}
  which is nondegenerate except for a $2$-dimensional radical in the fibre $\check{E}_P$.
  More precisely, $\check{\Sympl}_{2n, L, P}( S)$ is the following groupoid for each $k$-scheme $S$:
  \begin{itemize}
   \item Each object consists of a rank $2n$ vector bundle $\check{\stackE}$ on $C \times_k S$, a line bundle
    $\stackL$ on $C \times_k S$ locally in $S$ isomorphic to the pullback of $L$ from $C$, and an
    alternating form
$$
\check{b}: \check{\stackE} \otimes \check{\stackE} 
\longto  \stackL
$$
whose adjoint $b^{\#}:
    \check{\stackE} \longto  \check{\stackE}^{\dual} \otimes 
\stackL$ is injective with $\coker( b^{\#})$ locally,
    in $S$, isomorphic to the pullback of $k_P^2$ from $C$.
   \item Each morphism from $(\check{\stackE}, \check{b}: 
\check{\stackE} \otimes \check{\stackE} \longto  \stackL
    )$ to $(\check{\stackE}', \check{b}': \check{\stackE}' \otimes 
\check{\stackE}' \longto  \stackL')$ is a pair
    of vector bundle isomorphisms $\check{\stackE} \longto  
\check{\stackE}'$ and $\stackL \longto  \stackL'$ that
    intertwine $\check{b}$ and $\check{b}'$.
  \end{itemize}
  Over this stack $\check{\Sympl}_{2n, L, P}$, let
  \begin{equation} \label{Phicheck}
    \check{\Phi}: \Substack_{1, \O}( \Echeckuniv) \longto \check{\Sympl}_{2n, L, P}
  \end{equation}
 be the moduli stack of such bundles $(\check{E},\check{b})$ together 
with a trivial line subbundle
  of $\check{E}$ with the property that its fibre at $P$ is in the 
radical of $\check{E}_P$. More precisely, the stack
  $\Substack_{1, \O}( \Echeckuniv)$ is given by the following groupoid $\Substack_{1, \O}( \Echeckuniv)( S)$
  for each $k$-scheme $S$:
  \begin{itemize}
   \item Its objects consist of an object $(\check{\stackE}, \check{b}: \check{\stackE} \otimes \check{
    \stackE} \longto  \stackL)$ in $\check{\Sympl}_{2n, L, P}( 
S)$ and a line subbundle $\stackI \subseteq \check{
    \stackE}$, which is locally in $S$ isomorphic to $\O$, and for which 
the restriction
$$
\check{b}_P:
    \stackI_P \otimes \check{\stackE}_P \longto  \stackL_P
$$
to $\{P\} \times S$ vanishes identically.
   \item Morphisms from $(\check{\stackE}, \check{b}: \check{\stackE} 
\otimes \check{\stackE} \longto  \stackL,
    \stackI \subseteq \check{\stackE})$ to $(\check{\stackE}', \check{b}': \check{\stackE}' \otimes \check{
    \stackE}' \longto  \stackL', \stackI' \subseteq 
\check{\stackE}')$ consist of three vector bundle isomorphisms
$$
\check{\stackE} \longto  \check{\stackE}'\, , \quad \,\stackL 
\longto  
\stackL'
\,~\,\text{~and~}\, ~\, \stackI \longto  \stackI'
$$
that commute with $\check{b}$ and $\check{b}'$ and with the inclusions 
$\stackI \subseteq \check{\stackE}$ and $\stackI'
    \subseteq \check{\stackE}'$.
  \end{itemize}
  These objects have obvious scalar automorphisms; they turn $\check{\Sympl}_{2n, L, P}$ and $\Substack_{1, \O}( \Echeckuniv)$ into
  $\Gm$-stacks. Forgetting the subbundle $\stackI$ defines the $1$-morphism $\check{\Phi}$ in \eqref{Phicheck}; it is a morphism of weight $1$.

  We construct a diagram of stacks and $1$-morphisms over $k$
  \begin{equation}\label{phi.de}
 \xymatrix{
    \stackM \ar@{^{(}->}[r] \ar[d]_{\Phi|_{\stackM}} &
    \Substack_{1, \O}   ( \Euniv)      \ar[d]^{       \Phi}  
\ar[r]^-{\Pi_1} &
    \Substack_{1, \O}   ( \Echeckuniv) \ar[d]^{\check{\Phi}} \ar[r]^-{\Pi_2} &
    \Substack_{1, \O(P)}( \Etildeuniv) \ar[d]^{\widetilde{\Phi}}\\
    \stackU \ar@{}[r]|-{\displaystyle \subseteq} & \Sympl_{2n,L} & \check{\Sympl}_{2n,L,P} &
      \Sympl_{2n,L} \ar@{}[r]|-{\displaystyle \supseteq} & \widetilde{ \stackU}
  } \end{equation}
  as follows:
  \begin{itemize}
   \item $\Substack_{1, \O}( \Euniv)$ is the moduli stack introduced in Definition \ref{substack} that
    parameterises symplectic bundles $(E, b: E \otimes E \longto  
L)$ together with a trivial line subbundle
    $\O \subseteq E$.
   \item $\Substack_{1, \O(P)}( \Etildeuniv)$ is the moduli stack introduced in Definition \ref{substack}
    that parameterises symplectic bundles $( \widetilde{E}, \widetilde{b}: \widetilde{E} \otimes
    \widetilde{E} \longto  L)$ together with a subbundle $\O( P) 
\subseteq \widetilde{E}$.
   \item $\Phi := \Phi_{\O}$ and $\widetilde{\Phi} := \Phi_{\O( P)}$ are the forgetful $1$-morphisms of
    Definition \ref{substack}.
   \item $\Pi_1$ sends each triple $(E, b: E \otimes E \longto  
L, i: \O \subseteq E)$ to the kernel
    $\check{E} \subseteq E$ of the sheaf homomorphism
    \begin{equation*}
      E \longto[b^{\#}] \Hombdl( E, L) \longto[i^*] L \twoheadrightarrow L_P := L \otimes_{\O_C} k_P,
    \end{equation*}
    together with the restriction
$$
\check{b}: \check{E} \otimes  \check{E}\longto  L
$$ of $b$ and the same
    subbundle $\O \subseteq \check{E} \subseteq E$ (which is indeed contained in $\check{E}$ because $b$ is
    alternating).
   \item $\Pi_2$ sends each triple $(\check{E}, \check{b}: \check{E} 
\otimes \check{E} \longto  L, i: \O \subseteq
    \check{E})$ to the pushout $\O( P) \subseteq \widetilde{E}$ of $i$ along the sheaf monomorphism
    $\O \hookrightarrow \O( P)$, together with the alternating form $\widetilde{b}: \widetilde{E} \otimes
    \widetilde{E} \longto L$ that coincides with $\check{b}$ outside $P$.

    (A straightforward computation over the local ring $\O_{C, P}$ 
shows that there is precisely one such
    form $\widetilde{b}$ and that it is nondegenerate.)
   \item $\stackU \subseteq \Sympl_{2n, L}$ is the open substack of symplectic bundles $(E, b)$ with
    \begin{equation} \label{nonspecial}
      \dim \cohom^0( E) = n \qquad\text{and}\qquad \cohom^1( E) = 0
    \end{equation}
    for which the radical $\rad \cohom^0( E)$ of the pairing
    \begin{equation} \label{pair_sections}
      \cohom^0( E) \otimes \cohom^0( E) \xrightarrow{\eta_P \otimes \eta_P} E_P \otimes E_P \longto[b_P] L_P
    \end{equation}
    is $1$-dimensional, and any $0 \neq s \in \rad \cohom^0( E)$ is everywhere on $C$ nonzero.

    (Grothendieck's theory of cohomology and base change, \cite{ega3}, shows that \eqref{nonspecial} is an
    open condition. Assuming it, $\cohom^0( E)$ has odd dimension, so any alternating form on it is
    degenerate; hence $\dim \rad \cohom^0( E) = 1$ is then an open condition as well. This proves that
    $\stackU$ is indeed an open substack.)
   \item $\stackM \hookrightarrow \Phi^{-1}( \stackU)$ is the closed 
substack of triples
    $(E, b: E \otimes E \longto  L, \O \subseteq E)$ for which 
the restriction of \eqref{pair_sections} to
$$
\cohom^0(\O) \otimes \cohom^0( E) \subseteq \cohom^0(E) \otimes 
\cohom^0( E)
$$ vanishes.
   \item $\widetilde{\stackU} \subseteq \Sympl_{2n, L}$ is the open substack of symplectic bundles
    $(\widetilde{E}, \widetilde{b})$ with
    \begin{equation*}
      \dim \cohom^0( \widetilde{E}) = n \qquad\text{and}\qquad \cohom^1( \widetilde{E}) = 0
    \end{equation*}
    for which the evaluation map $\widetilde{\eta}_P: \cohom^0( 
\widetilde{E}) \longto  \widetilde{E}_P$ has rank
    $\geq n-1$.
  \end{itemize}
  By construction, $\Phi$ restricts to a $1$-isomorphism 
$\Phi|_{\stackM}: \stackM \longto  \stackU$; its inverse
  endows each symplectic bundle $(E, b)$ with the image $\O \subseteq E$ of any $0 \neq s \in \rad \cohom^0(
  E)$. Thus it suffices to show that $\stackM$ is non-empty and rational as a $\Gm$-stack.
  \begin{lemma} \label{non-empty}
    Even the intersection $\stackM \cap \Pi_1^{-1} \Pi_2^{-1} \widetilde{ \Phi}^{-1}( \widetilde{\stackU})$ is non-empty.
  \end{lemma}
  \begin{proof}
    We may assume without loss of generality that $k = \bar{k}$ is algebraically closed.
    Let $F$ be a general vector bundle extension
    \begin{equation*}
      0 \longto \O \longto F \longto L \longto 0,
    \end{equation*}
    and let the vector bundle extension
    \begin{equation*}
      0 \longto \O( P) \longto \widetilde{F} \longto L( -P) \longto 0
    \end{equation*}
    be the image of $[F]$ under the natural surjection
$$\Ext^1( L, \O) \twoheadrightarrow \Ext^1( L( -P),
    \O( P))\, .
$$
According to Riemann-Roch and our assumptions on $L$, we have
    \begin{equation*}
      \dim \cohom^0( L) = g = \dim \cohom^1( \O) \quad\text{and}\quad
      \dim \cohom^0( L( -P)) = g-1 = \dim \cohom^1( \O( P));
    \end{equation*}
    using Lemma \ref{nonspecial_extension}, it follows that
    \begin{equation*}
      \cohom^0( F) \cong k \cong \cohom^0( \widetilde{F}) \quad\text{and}\quad
      \cohom^1( F)   =   0   =   \cohom^1( \widetilde{F}).
    \end{equation*}
    Now let $(E', b': E' \otimes E' \longto  L)$ be a symplectic 
bundle of rank $2n-2$ with the properties given
 in Lemma \ref{generic}, and let $(E, b)$ be the fibrewise 
orthogonal direct sum of $E'$ and $F$. Then the tuple
    \begin{equation*}
      (E, \quad b, \quad \O \cong 0 \oplus \O \subseteq E' \oplus F = E)
    \end{equation*}
    defines a point in $\stackM$; its image under $\widetilde{ \Phi} \circ \Pi_2 \circ \Pi_1$ is by
    construction the symplectic bundle $\widetilde{E} := E' \perp \widetilde{F}$ and thus contained in
    $\widetilde{ \stackU}$.
  \end{proof}
  
According to Proposition \ref{unirat}(ii), the $\Gm$-stack $\Substack_{1, \O(P)}( \Etildeuniv)$ is rational.
  To deduce from that the required rationality of $\stackM$, we study the fibres of the above $1$-morphisms
  $\Pi_1$ and $\Pi_2$ via the strictly commutative diagram
  \begin{equation*} \xymatrix{
    \stackM \hookrightarrow \Substack_{1, \O}( \Euniv) \ar[r]^-{\Gamma_1} \ar[d]_{\Pi_1} &
    \P( \rad \Echeckuniv_P) \ar[dl] &\\
    \Substack_{1, \O}( \Echeckuniv) \ar[r]^-{\Gamma_2} \ar[dr]_{\Pi_2} & \PP( \Etildeuniv_P) \ar[d] &
    \check{\stackM} = \PP( \frac{\Etildeuniv_P}{\im \cohom^0( \Etildeuniv)}) \ar@{_{(}->}[l] \ar[d] &\\
    & \Substack_{1, \O( P)}( \Etildeuniv) &
    \widetilde{ \Phi}^{-1}( \widetilde{ \stackU}) \ar@{}[l]|-{\displaystyle \supseteq}
  } \end{equation*}
  constructed as follows:
  \begin{itemize}
   \item $\rad \Echeckuniv_P$ is the vector bundle of rank $2$ and weight $1$ on $\Substack_{1, \O}(
    \Echeckuniv)$ whose fibre over any triple $(\check{E}, \check{b}: 
\check{E} \otimes \check{E} \longto  L, \O
    \subseteq \check{E})$ is the radical $\rad \check{E}_P$ of the alternating pairing
    $\check{b}_P: \check{E}_P \otimes \check{E}_P \longto  L_P$ 
on the fibre $\check{E}_P$.
   \item $\Gamma_1$ sends a triple $(E, b: E \otimes E \longto  
L, \O \subseteq E)$ to its image $\check{E}
    \subseteq E$ under $\Pi_1$ and the image of the $k$-linear map 
$E(-P)_P \longto  \rad \check{E}_P \subseteq
    \check{E}_P$ induced by the inclusion $E( -P) \subseteq \check{E}$ as subsheaves of $E$.
   \item $\Etildeuniv_P$ is the vector bundle of rank $2n$ and weight $1$ on $\Substack_{1, \O( P)}(
    \Etildeuniv)$ whose fibre over any triple $(\widetilde{E}, 
\widetilde{b}, \O( P) \subseteq \widetilde{E})$
    is the fibre $\widetilde{E}_P$ of $\widetilde{E}$ at $P$.
   \item $\Gamma_2$ sends a triple $(\check{E}, \check{b}, \O \subseteq 
\check{E})$ to its image $\check{E}
    \subseteq \widetilde{E}$ under $\Pi_2$ and the image of the induced 
$k$-linear map $\check{E}_P \longto 
    \widetilde{E}_P$.
   \item $\check{\stackM} \hookrightarrow \PP( \Etildeuniv_P)|_{\widetilde{\Phi}^{-1}( \widetilde{\stackU})}$
    is the closed substack of those hyperplanes in the fibres $\widetilde{E}_P$ which contain
    the values at $P$ of all section $s \in \cohom^0( \widetilde{E})$.
  \end{itemize}
  The vector spaces $\cohom^0( \widetilde{E})$ are the fibres of a rank $n$ vector bundle $\cohom^0(
  \Etildeuniv)$ on $\widetilde{\Phi}^{-1}( \widetilde{\stackU})$ by cohomology and base change \cite{ega3}. 
  The evaluations $\widetilde{\eta}_P: \cohom^0( \widetilde{E}) 
\longto  \widetilde{E}_P$ at $P$ define a morphism
  of vector bundles on $\widetilde{\Phi}^{-1}( \widetilde{\stackU})$
  \begin{equation*}
    \widetilde{\eta}_P^{\univ}: \cohom^0( \Etildeuniv) \longto \Etildeuniv_P.
  \end{equation*}
  Every section $s \in \cohom^0( \O( P))$ vanishes in the fibre $\O( P)_P$. If $\widetilde{E}$ admits a
  subbundle isomorphic to $\O( P)$, then the rank of $\widetilde{\eta}_P$ can thus be at most $n-1$. This
  and the definition of $\widetilde{\stackU}$ show that $\widetilde{\eta}_P^{ \univ}$ has constant rank $n-1$
  on $\widetilde{\Phi}^{-1}( \widetilde{\stackU})$; thus its cokernel
  \begin{equation*}
    \coker( \widetilde{\eta}_P^{\univ}) = \frac{\Etildeuniv_P}{\im \cohom^0( \Etildeuniv)}   
  \end{equation*}
  is a vector bundle of weight $1$ and rank $n+1$ on $\widetilde{\Phi}^{-1}( \widetilde{\stackU})$. The
  associated projective subbundle of $\PP( \Etildeuniv_P)|_{\widetilde{\Phi}^{-1}( \widetilde{\stackU})}$
  is by construction the above closed substack $\check{\stackM}$.

  Now Lemma \ref{non-empty} implies in particular that $\widetilde{\Phi}^{-1}( \widetilde{\stackU})$ is
  non-empty. It follows that the projective bundle $\check{\stackM}$ over it is also non-empty, algebraic,
  locally of finite type over $k$, smooth, irreducible, and rational as a $\Gm$-stack.

  It is easy to check that $\Gamma_1$ and $\Gamma_2$ are open immersions. More precisely:
  \begin{itemize}
   \item $\Gamma_1$ is a $1$-isomorphism onto the open substack of
tuples $(\check{E}, \check{b}, \O
    \subseteq \check{E})$ together with a line in $\rad \check{E}_P$ different from the fibre of $\O
    \subseteq \check{E}$. Its inverse can be described as follows:

    Using the canonical isomorphism $\P \check{E}_P \cong \P \check{E}(P)_P$, such a line in $\check{E}_P$
    yields a line in $\check{E}(P)_P$. Let $E \subseteq \check{E}(P)$ be its inverse image under the sheaf
    surjection $\check{E}( P) \twoheadrightarrow \check{E}( P)_P$. Then it is easy to see that $\O \subseteq
    \check{E}$ is in fact a subbundle of $E$. A straightforward 
computation over the local ring $\O_{C, P}$
    shows that there is precisely one alternating form $b$ on $E$ equal to $\check{b}$ outside $P$, and that
    $b$ is nondegenerate. Now the inverse of $\Gamma_1$ sends $(\check{E}, \check{b}, \O \subseteq \check{E})
    $ and the line in question to $(E, b, \O \subseteq E)$.

   \item $\Gamma_2$ is a $1$-isomorphism onto the open substack of
triples $(\widetilde{E}, \widetilde{b},
    \O( P) \subseteq \widetilde{E})$ together with a hyperplane in $\widetilde{E}_P$ which does not contain
    the fibre of the subbundle $\O( P) \subseteq \widetilde{E}$. Its inverse sends such a hyperplane to its
    inverse image $\check{E} \subseteq \widetilde{E}$ under the sheaf surjection $\widetilde{E}
    \twoheadrightarrow \widetilde{E}_P$, equipped with the restricted form
    $\check{b} := \widetilde{b}|_{\check{E} \otimes \check{E}}$ and the subbundle
    $\O \cong \O(P) \cap \check{E} \subseteq \check{E}$.
  \end{itemize}
  This implies in particular that the stack $\Substack_{1, \O}( \Echeckuniv)$ is algebraic, locally of finite
  type over $k$, smooth and irreducible.

  Inside the open substack
  \begin{equation} \label{open}
    \Phi^{-1}( \stackU) \cap \Pi_1^{-1} \Pi_2^{-1} 
\widetilde{\Phi}^{-1}( \widetilde{\stackU})
      \subseteq \Substack_{1, \O}( \Euniv),
  \end{equation}
(see \eqref{phi.de}) we have the closed substack
  \begin{equation} \label{M_int}
    \stackM \cap \Pi_1^{-1} \Pi_2^{-1} \widetilde{\Phi}^{-1}( 
\widetilde{\stackU})
  \end{equation}
  which has been defined by the vanishing of
$$
b_P \circ (\eta_P \otimes \eta_P): \cohom^0(\O) \otimes
  \cohom^0(E) \longto  L_P\, .
$$
We also have the closed substack
  \begin{equation} \label{Mcheck_int}
    \Phi^{-1}( \stackU) \cap \Pi_1^{-1} \Gamma_2^{-1}( \check{ \stackM})
 \subseteq \Phi^{-1}( \stackU) \cap \Pi_1^{-1} \Pi_2^{-1} 
\widetilde{\Phi}^{-1}( \widetilde{\stackU})
\end{equation}
of all triples $(E, b, \O \subseteq E)$ such that their images $E 
\supseteq 
\check{E} \subseteq \widetilde{E}$
  under $\Pi_1$ and $\Pi_2$ satisfy the following closed condition:
  \begin{equation*}
    \widetilde{\eta}_P( s) \in \widetilde{E}_P \text{ lies in the image 
of } \check{E}_P \longto  \widetilde{E}_P
    \text{ for all } s \in \cohom^0( \widetilde{E}).
  \end{equation*}
  Now these two closed conditions are equivalent: Since $\dim \cohom^0(E) =n= \dim \cohom^0( \widetilde{E})$
  everywhere on the open substack \eqref{open}, both closed conditions are
equivalent to the condition that $\dim
  \cohom^0( \check{ E}) = n$. Hence the closed substacks \eqref{M_int} and \eqref{Mcheck_int} coincide.

  In particular, \eqref{Mcheck_int} is non-empty, because \eqref{M_int} is so according to Lemma
  \ref{non-empty}. But $\Gamma_1$ yields an open immersion of \eqref{Mcheck_int} into the restriction of the
  projective bundle $\P( \rad \Echeckuniv_P)$ to the rational $\Gm$-stack $\Gamma_2^{-1}( \check{\stackM})$.
  It follows that the $\Gm$-stack \eqref{Mcheck_int} is rational as well; hence the same holds for
  \eqref{M_int} and consequently also for $\stackM$.
\end{proof}

\begin{cor}
  Under the hypotheses of Theorem \ref{theorem1}, the coarse moduli scheme $\Symcoarse_{2n, L}$ is rational.
\end{cor}

\begin{rem} \label{automatic}
  Suppose that $C$ admits a rational point $P \in C( k)$ and a line bundle $\xi$ of degree $0$ on $C$ with
  $\xi^{\otimes 2} \not\cong \O_C$. Then there is, for every given line bundle $L$ of odd degree on
  $C$, a line bundle $I$ with the properties required in Theorem \ref{theorem1}. Indeed, we may take either $I := \O_C( dP)$
  or $I := \xi( dP)$, where $d := (2g - 1 - \deg L)/2 \in \integers$.
\end{rem}

\begin{rem}
  One can also vary the line bundle $L$, keeping only its degree fixed. More precisely, let
  \begin{equation*}
    \Sympl_{2n, d}
  \end{equation*}
  denote the moduli stack of all tuples $(E, b: E \otimes E 
\longto  L)$ in which $E$ is a vector bundle of rank
  $2n$ on $C$, $L$ is a line bundle of degree $d$ on $C$, and $b$ is a nowhere degenerate symplectic form.
  Let $\Symcoarse_{2n, d}$ be the corresponding coarse moduli scheme of Ramanathan-stable symplectic bundles $( E, b: E \otimes E \longto L)$,
  as constructed in \cite{ramanathan} and \cite{GLSS}.
  Concerning their birational type, we have the following:
\end{rem}
\begin{cor} \label{var_det}
  Suppose that $n$ and $d$ are both odd, and that the curve $C$ over $k$ has a rational point $P \in C( k)$.
  Then the stack $\Sympl_{2n, d}$ is birational to
  \begin{equation*}
    \BGm \times \Aa^s \times \Pic^d( C) \quad  \text{ for some $s$.} 
  \end{equation*}
\end{cor}
\begin{proof}
  Forgetting $E$ and $b$ defines a canonical $1$-morphism
  \begin{equation*}
    \Sympl_{2n, d} \longto \Pic^d( C)
  \end{equation*}
  to the Picard \emph{scheme} $\Pic^d( C)$. The fibres of this morphism are the moduli stacks $\Sympl_{2n,L}$
  studied above. In particular, its generic fibre coincides with the moduli stack $\Sympl_{2n, L^{\generic}}$
  over the function field $K := k( \Pic^d(C))$, where the line bundle $L^{\generic}$ on $C_K := C \times_k K$
  is the generic fibre of a Poincar\'{e} family on $C \times_k \Pic^d( C)$. Now $C_K$ has a rational point
  because $C$ has; due to Remark \ref{automatic}, it thus suffices to construct a line bundle $\xi$ of degree
  $0$ on $C_K$ with $\xi^{\otimes 2} \not\cong \O_{C_K}$.

  We take for $\xi$ the generic fibre of a Poincar\'{e} family on $C \times_k \Pic^0( C)$.
  Then $\xi^{\otimes 2} \not\cong \O_{C_K}$ holds indeed, because the endomorphism of the abelian variety $\Pic^0( C)$
  that sends each line bundle on $C$ to its square is nonconstant.
\end{proof}

\begin{rem}
  Corollary \ref{var_det} automatically implies the following: $\Symcoarse_{2n, d}$ is birational to $\Aa^s \times \Pic^d( C)$
  for some $s$ if $n \cdot d$ is odd and $C( k)$ is non-empty.
\end{rem}


\begin{thebibliography}{11}

\bibitem{GLSS}
T.L. Gomez, A. Langer, A.H.W. Schmitt, I. Sols,
Moduli Spaces for Principal Bundles in Arbitrary Characteristic,
Adv. Math. (2008), doi:10.1016/j.aim.2008.05.015.

\bibitem{ega3}
A.~Grothendieck,
EGA III: \'{E}tude cohomologique des faisceaux coh\'erents,
Publ. Math. Inst. Hautes \'Etudes Sci. 11 (1961);
Publ. Math. Inst. Hautes \'Etudes Sci. 17 (1963).

\bibitem{sga4}
A.~Grothendieck, et~al.,
SGA IV: Th\'eorie des topos et cohomologie \'etale des
  sch\'emas, Lecture Notes in Math., vols. 269, 270, 305,
Springer-Verlag, Berlin, 1972/73.

\bibitem{par}
N. Hoffmann,
Rationality and Poincar\'e families for vector bundles with extra
  structure on a curve, Int. Math. Res. Not. 2007
(2007), article ID rnm010.

\bibitem{KS}  A. King, A. Schofield, Rationality of moduli of vector 
bundles on curves, Indag. Math. (N.S.) 10 (1999) 519--535.

\bibitem{Ne} P.E. Newstead, Rationality of moduli spaces of stable 
bundles, Math. Ann. 215 (1975) 251--268.

\bibitem{ramanathan}
A. Ramanathan, Moduli for principal bundles over algebraic curves,
Proc. Indian Acad. Sci. Math. Sci. 106 (3) (1996) 301--328,
Proc. Indian Acad. Sci. Math. Sci. 106 (4) (1996) 421--449.
\end{thebibliography}
\end{document}